\documentclass[oneside,10pt]{amsart}
\usepackage{amsmath,amsfonts, latexsym,amssymb}
\usepackage[active]{srcltx}

\theoremstyle{definition}

\numberwithin{equation}{section}


\begin{document}

{\renewcommand{\thefootnote}{}\footnote{{\bf 2010 Mathematics Subject 
Classification:} 46E10, 46J15, 46J20, 30H50, 30H15.\\   
\noindent{\bf Keywords and Phrases:}  class $M^p$ $(1<p<\infty)$,
     Privalov class $N^p$ $(1<p<\infty)$, 
 $F$-algebra, prime  ideal, principal ideal, maximal ideal.}

\title[A characterization of some prime ideals in certain  
 $F$-algebras \ldots]
{A characterization of some  prime ideals in certain  
 $F$-algebras of holomorphic functions}

\author{Romeo Me\v{s}trovi\'{c}}
\address{Maritime Faculty Kotor, University of Montenegro, Dobrota 36,
85330 Kotor, Montenegro} \email{romeo@ucg.ac.me}

\address{Maritime Faculty Kotor, University of Montenegro, 
Dobrota 36, 85330 Kotor, Montenegro} 

\email{romeo@ucg.ac.me}

\begin{abstract} The class $M^p$ $(1<p<\infty)$
 consists of all holomorphic functions $f$ on 
the open unit disk $\Bbb D$ for  which 
 $$
\int_0^{2\pi}\left(\log^+Mf(\theta)\right)^p\,\frac{d\theta}{2\pi}<\infty,
$$
where $Mf(\theta)=\sup_{0\leqslant r<1}\big\lvert f\big(re^{i\theta}\big)\big
\rvert$.
The class $M^p$ equipped with the topology given by the metric 
$\rho_p$ defined by            
 $$\rho_p(f,g)=\Vert f-g\Vert_p=\left(\int_0^{2\pi}\log^p(1+M(f-g)(\theta))\,
\frac{d\theta}{2\pi}\right)^{1/p}\quad (f,g\in M^p)$$
becomes an $F$-algebra.

In this paper, we consider the ideal structure of the classes 
$M^p$ $(1<p<\infty)$. Our main result   gives a complete characterization 
of prime ideals in  $M^p$ which are not  dense subsets of $M^p$. 
As a consequence, we obtaiin a related Mochizuki's  result 
concerning the Privalov classes  $N^p$ $(1<p<\infty)$.
  \end{abstract}

\maketitle

\vspace{-2mm}

\section{Introduction and Preliminaries}

Let $\Bbb D$ denote the open unit disk 
in the complex plane  and let $\Bbb T$ denote the
boundary of $\Bbb D$.  Let $L^q(\Bbb T)$ $(0<q\le \infty)$ be the
familiar Lebesgue space on the unit circle $\Bbb T$.

Following  Kim \cite{k1,k2}, 
the  class $M$ consists of all holomorphic functions $f$
on $\Bbb D$ for which
   $$
\int_0^{2\pi}\log^+Mf(\theta)\frac{d\theta}{2\pi}<\infty,\eqno(1)
     $$
where $\log^+|a|=\max\{\log |a|,0\}$ and
   $$
Mf(\theta)=\sup_{0\leqslant r<1}\big\lvert f\big(re^{i\theta}\big)\big\rvert
\eqno(2)
  $$
is the  maximal radial function of $f$.

The study of the class $M$ has been 
well established in \cite{k1}, \cite{k2},
\cite{gz} and M. Nawrocky \cite{n}.  
Kim \cite[Theorems 3.1 and 6.1]{k2} showed that the space $M$ with 
the topology given by the metric $\rho$ defined by
    $$
\rho(f,g)=\int_0^{2\pi}\log(1+M(f-g)(\theta))\,\frac{d\theta}{2\pi},
\quad f,g\in M,\eqno (3)
    $$    
becomes an $F$-algebra.

 Recall that the  Smirnov class $N^+$  is the set
of all functions  $f$ holomorphic on $\Bbb D$ such that
   $$
\lim_{r\rightarrow 1}\int_0^{2\pi}\log^+\vert f(re^{i\theta})\vert
\,\frac{d\theta}{2\pi}=\int_0^{2\pi}\log^+\vert f^*(e^{i\theta})\vert\,
\frac{d\theta}{2\pi}<+\infty,\eqno(4)
   $$
where $f^*$ is the boundary function of $f$ on $\Bbb T$; that is,
 $$
f^*(e^{i\theta})=\lim_{r\rightarrow 1-}f(re^{i\theta})\eqno(5)
  $$ 
is the  radial  limit of $f$  which exists for almost 
every $e^{i\theta}\in\Bbb T$.

The classical  Hardy space $H^q$ $(0<q\le\infty)$
 consists of all functions $f$ holomorphic on $\Bbb D$, which satisfy
   $$
\sup_{0\le r<1}\int_0^{2\pi}\big\lvert f\big(re^{i\theta}\big)\big\rvert^q
\frac{d\theta}{2\pi}<\infty\eqno(6)
   $$
if $0<q<\infty$, and which are bounded when $q=\infty$:
   $$
\sup_{z\in D}\lvert f(z)\rvert<\infty.\eqno(7)
    $$

Although the class $M$ is essentially smaller than the class $N^+$,
 it was showed in \cite{n}  
 that the class $M$ and the Smirnov class $N^+$ have the same corresponding 
locally convex structure which was already established for $N^+$ by 
N. Yanagihara for the Smirnov class in   \cite{y2} and \cite{y1}
(also see \cite{mc}).

The  Privalov class $N^p$ $(1<p<\infty)$ is defined 
as the set of all holomorphic functions $f$ $f$ on $\Bbb D$ such that 
   $$
\sup_{0\le  r<1}\int_0^{2\pi}(\log^+\vert f(re^{i\theta})\vert)^p\,
\frac{d\theta}{2\pi}<+\infty\eqno(9)
   $$ 
holds. These classes were firstly considered  by Privalov in 
\cite[p. 93]{p}, where $N^p$ is denoted as $A_q$. 

Notice that for $p=1$, the condition (9) defines the 
 Nevanlinna class $N$ of holomorphic functions in $\Bbb D$.

It is known (see \cite{mo, mp1, i1, ck})  that the following inclusion 
relations hold: 
   $$
N^r\subset N^p\;(r>p),\quad\bigcup_{q>0}H^q\subset
\bigcap_{p>1}N^p,\quad{\rm and}\quad\bigcup_{p>1}N^p\subset M\subset N^+
\subset N,\eqno(10)
   $$
where the above containment relations are proper.

The study of the spaces $N^p$ $(1<p<\infty)$ was continued in 1977
by M. Stoll  \cite{s} (with the notation
 $(\log^+H)^\alpha$ in \cite{s}).  Further, the topological and functional 
properties of these spaces  have been  studied by several authors
(see \cite{mo}--\cite{msu}).

It is well known (see, e.g., \cite[p. 26]{d}) that a function $f\in N$
belongs to $N^{+}$ if and only if 
  $$
f(z)=e^{i\gamma}B(z)S(z)F(z),z\in\Bbb D,\eqno(11)
  $$
where $\gamma$ is a real constant,  $B$ is a Blaschke product with respect 
to zeros of $f(z)$, $S$ is a singular inner function and $F$ is an outer 
function for the class $N$, i.e.,
   $$
B(z)=z^m\prod_{n}\frac{|a_n|}{a_n}\cdot\frac{a_n-z}{1-\overline{a_n}z},
\quad z\in\Bbb D.\eqno(12)
  $$
where $m$ is a nonnegative integer and $\sum_{n}(1-|a_n|)<\infty$,
  $$
S(z)=\exp\left(-\int_0^{2\pi}\frac{e^{it}+z}{e^{it}-z}\,
d\mu_k(t)\right),\quad z\in\Bbb D,\eqno(13)
  $$
with positive singular measure $d\mu$,  and 
    $$
F(z)=\omega\exp\left(\frac1{2\pi}\int_0^{2\pi}
\frac{e^{it}+z}{e^{it}-z}\log\big\lvert f^*\big(e^{it}\big)\big\rvert\,
dt\right),\quad z\in\Bbb D,\eqno(14)
  $$
where $\log\vert F^*\vert\in L^1(\Bbb T)$.

Recall that a function $I$ of the form 
$$
I=e^{i\gamma}BS\eqno(15)
 $$ 
is called an inner function. It is well known  
(see, e.g., \cite[p.24]{d}) that  $I$ is a bounded holomorphic function
on $\Bbb D$ such that $|I(z)|\le 1$ for each $z\in\Bbb D$,   
$|I^*(e^{i\theta})|=1$ and $|B^*(e^{i\theta})|=1$  for almost every $e^{i\theta}\in\Bbb T$.

Privalov \cite{p} showed the following 
canonical factorization theorem for the classes $N^p$ ($1<p<\infty$). 

\vspace{2mm}
\noindent{\bf Theorem 1} (\cite[pp. 98--100]{p}; also see 
\cite{e2}). {\it Let $p>1$ be a fixed real number. 
A function $f\in N^p$ $(f\not\equiv 0)$ has a unique factorization of the form
    $$
f(z)=B(z)S(z)F(z),\eqno(16)
   $$
where $B$ is the Blaschke product with respect to zeros of $f(z)$,
$S(z)$ is a singular inner function and $F(z)$ is an outer function
such that $\log^+\vert F^*\vert\in L^p(\Bbb T)$.
Conversely, every such product $BSF$ belongs to $N^p$.}

\vspace{2mm}
M. Stoll  proved the following result.

\vspace{2mm}
\noindent{\bf Theorem 2} (\cite[Theorem 4.2]{s}). 
{\it For each $p>1$ the Privalov space $N^p$ with the topology given by the metric 
$d_p$ defined by            
   $$
d_p(f,g)=\Big(\int_0^{2\pi}\big(\log(1+
\vert f^*(e^{i\theta})-g^*(e^{i\theta})\vert)\big)^p\,\frac{d\theta}
{2\pi}\Big)^{1/p},\quad f,g\in N^p,\eqno (17) 
   $$  
becomes an $F$-algebra, that is, $N^p$ is an  $F$-space  in which  
multiplication is continuous.}

\vspace{2mm}

Recall that the function $d_1=d$ defined on the Smirnov class $N^+$
by (17) with $p=1$ induces the metric topology on $N^+$.
N. Yanagihara \cite{y1} showed that 
under this topology, $N^+$ is an $F$-space.

Motivated by the mentioned investigations of the classes
$M$ and $N^{+}$, and the fact that the classes
$N^p$ $(1<p<\infty)$ are generalizations 
of the Smirnov class $N^+$, in  \cite{m1}, \cite{m2} and  \cite{m3} the author of this paper
investigated the classes  $M^p$ $(1<p<\infty)$
as generalizations of the class $M$. Accordingly, the class $M^p$
$(1<p<\infty)$ consists of all holomorphic functions $f$ on 
$\Bbb D$ for  which 
    $$
\int_0^{2\pi}\left(\log^+Mf(\theta)\right)^p\,\frac{d\theta}{2\pi}<\infty.
\eqno (18)
    $$
Obviously, $\bigcup_{p>1}M^p\subset M$ holds.
Following \cite{m3}, by analogy with the space $M$, the space $M^p$ can 
be equipped with the topology 
induced by the metric $\rho_p$ defined as            
 $$
\rho_p(f,g)=\left(\int_0^{2\pi}\log^p(1+M(f-g)(\theta))\,
\frac{d\theta}{2\pi}\right)^{1/p},\eqno(19)
 $$
with $f,g\in M^p$. 

For our purposes, we will also need the following two results.

\vspace{2mm}

\noindent{\bf Theorem 3} (\cite[Corollary 18]{m3}. {\it For each $p>1$,
 $M^p$ is an $F$-algebra with respect to the metric topology induced 
by the metric $\rho_p$ defined by $(19)$.}

\vspace{2mm}

\noindent{\bf Theorem 4} (\cite[Theorem 16]{m3}. {\it For each $p>1$,
the classes $M^p$ and $N^p$ coincide, and the metric spaces 
$(M^p,\rho_p)$ and $(N^p,d_p)$ have the same topological structure.}

\vspace{2mm}

Since the space $M^p$ $(1<p<\infty)$ is an algebra, it can be 
also  considered as a ring with respect  to the usual  
ring's operations addition and multiplication. Notice that 
these two operations are continuous in $M^p$  because $M^p$ is an $F$-algebra. 

Motivated by several results on the ideal structure of some spaces 
of holomorphic functions given in \cite{k1}, \cite{mo}, \cite{ma} and 
\cite{be}-\cite{m17},
related investigations for the spaces $N^p$ $(1<p<\infty)$
and their Fr\'{e}chet envelopes were given in 
\cite{mo},  \cite{me4}, \cite{ma},    \cite{mp5}, \cite{mp3} and \cite{m5}.
Note that a survey of these results was given in \cite{mp6}.
A complete characterization of principal  ideals 
in $N^p$ which are dense in $N^p$ was given by Mochizuki \cite[Theorem 3]{mo}. 
This result was used in \cite{mp5} to obtain the
 $N^p$-analogue of the famous  Beurling's theorem 
for the Hardy spaces $H^q$ $(0<q<\infty)$ \cite{be}.
Moreover, it was proved in \cite[Theorem B]{me4}
that $N^p$ $(1<p<\infty)$ is a ring of Nevanlinna--Smirnov type  
in the sense of Mortini \cite{Mor}. The structure of closed weakly 
dense ideals in $N^p$ was established in \cite{mp3}.
The structure  of maximal ideals in $N^p$
was studied in \cite{m5}, where related results are similar to these  
obtained by Roberts and Stoll \cite{RS} for the Smirnov class $N^+$.

Let $R$ be a commutative ring.
An ideal  ${\mathcal I}$ in  $R$ is called principal if there is an 
element $a$ of $R$ such that ${\mathcal I}=aR:=\{ab: b\in R \}$.
 An ideal  ${\mathcal I}$ in  $R$ is called maximal
if  ${\mathcal I}\not= R$ and no proper ideal of $R$ properly contains
${\mathcal I}$. An ideal  ${\mathcal I}$ in $R$ is called a prime ideal if for
any $a,b\in R$ if $a,b\in {\mathcal I}$ then either $a\in {\mathcal I}$
or $b\in {\mathcal I}$.     
Our goal in this paper is to characterize prime ideals in the 
algebras $M^p$ $(1<p<\infty)$ which are not  dense subsets of $M^p$
  (Theorem 5 in the next section). 
As an application  
of this result and Theorem 4, we obtain a related  result of Mochizuki
\cite{mo} concerning a characterization 
of closed prime ideals in the algebras  $N^p$ $(1<p<\infty)$
which are not dense in $N^p$  (Corollary 6).

\section{The main result and its proof}

Let $p>1$ be any fixed real number. For $\lambda\in \Bbb D$, we define    
  $$
{\mathcal M}_\lambda =\{f\in M^p:f(\lambda)=0\}.\eqno(20)
  $$
It was proved in \cite[Proposition 2.2]{msu} that for each $\lambda\in\Bbb D$,
a set ${\mathcal M}_\lambda$
defined  by (20) with $N^p$ instead of $M^p$,  is a 
closed maximal ideal in $N^p$ and that ${\mathcal M}_\lambda = (z-\lambda)$ 
\cite[Proof of Theorem 2.3]{msu}. 

The main result of this paper (Theorem 5) is the $M^p$-analogue of Theorem 4 
in \cite{k1} concerning a characterization of prime ideals in  $M$ 
which are not dense in $M$.

\vspace{2mm}

\noindent{\bf Theorem 5.} {\it   Let  ${\mathcal M}$ be a nonzero 
prime ideal in   $M^p$ $(1<p<\infty)$ which is not a dense subset of $M^p$. 
Then ${\mathcal M}={\mathcal M}_{\lambda}$, for some 
$\lambda\in \Bbb D$, where ${\mathcal M}_\lambda$ is defined by $(20)$.}
\vspace{2mm}

As an immediate consequence of Theorems 4 and  5, we immediatelyly obtain
the following $N^p$-analogue of Theorem 5 established by N. Mochizuki \cite{mo}.  

\vspace{2mm}

\noindent {\bf Corollary 6} (\cite[Theorem 35]{mo}). {\it  Let  ${\mathcal M}$
 be a nonzero 
prime ideal in  $N^p$ $(1<p<\infty)$ which is not a dense subset of $N^p$. 
Then ${\mathcal M}={\mathcal M}_{\lambda}$, for some 
$\lambda\in \Bbb D$, where ${\mathcal M}_\lambda$ is defined by $(20)$ with
 $N^p$ instead of $M^p$.}
\vspace{2mm}

Notice that the $N^+$-analogue of Theorem 5 was proved by J. Roberts and 
M. Stoll \cite[Theorem 1]{RS}. 

Proof of Theorem 5 is similar to those  of Theorem 4 in 
\cite{k1}, and it is  based on the following five lemmas.

\vspace{2mm}

\noindent{\bf Lemma 7.} {\it  Let  ${\mathcal M}$
be a nonzero  ideal od $M^p$. Then ${\mathcal M}$ contains a bounded 
holomorphic function which is not identically zero.} 

\vspace{2mm}

\noindent{\it Proof.}
Let  $f\in {\mathcal M}$ and  suppose that $f\not\equiv 0$. 
By Theorems 1 and 4,   $f$ can be factorized as
      $$
f(z)=B(z)S(z)F(z),\quad z\in \Bbb D,\eqno(21)
      $$
where $B$ is  the Blaschke product with respect to zeros 
of  $f$, $S$ is the  singular inner function  and 
      $$
F(z)=\exp\left(\frac{1}{2\pi}\int_0^{2\pi}\frac{e^{it}+z}{e^{it}-z}
\log |f^*(e^{it})|\,dt\right),\quad z\in \Bbb D,\eqno(22)
     $$
is the outer function. Define the function $g$ as
     $$
g(z)=\exp\left(-\frac{1}{2\pi}\int_0^{2\pi}\frac{e^{it}+z}{e^{it}-z}
\log^+|f^*(e^{it})|\,dt\right),\quad z\in \Bbb D,\eqno(23)    
     $$
Then $g$ is a bounded holomorphic function on $\Bbb D$, and thus, 
 $g\in N^p$. Since  ${\mathcal M}$ is an ideal of $M^p$,
it follows that $fg\in {\mathcal M}$. On the other hand, we have 
    $$
f(z)g(z)=B(z)S(z)\exp\left(\frac{1}{2\pi}\int_O^{2\pi}\frac{e^{it}+z}{e^{it}-z}
\log^-|f^*(e^{it})|\,dt\right),\quad z\in \Bbb D,\eqno(24)
    $$
where $\log^-|a|=\log^+|a|-\log|a|$. From the factorization  (24)
it follows that  $|f(z)g(z)|<1$ for each $z\in \Bbb D$, 
and therefore, $fg\in H^\infty$.  Thus, the function $fg$
has the desired  property.  \hfill $\square$ 
\vspace{2mm}

\noindent{\bf Lemma 8.} {\it  Let $F\in M^p$ be a function such that
 $F(z)\not=0$ for each  $z\in \Bbb D$. Then there exists 
a sequence  $\{F_n\}_{n=1}^{\infty}$ of functions in $M^p$ 
such that  $\left(F_n\right)^n=F$ for all positive integers  $n$
and $F_n\to 1$ in the space  $M^p$ as $n\to\infty$.}

\vspace{2mm}

\noindent{\it Proof.} 
Since by the assumption, $F(z)\not=0$ for each  $z\in \Bbb D$,
there exist a real-valued function   
$\rho(z)$ on $\Bbb D$ and a positive  
continuous  function  $\phi(z)$ on $\Bbb D$ for which
     $$
   F(z)=\rho(z)e^{i\phi(z)},\quad z\in \Bbb D.\eqno(25)
     $$
Now define the  sequence  $\{F_n\}_{n=1}^{\infty}$ of functions as 
     $$ 
F_n(z)=\rho(z)^{1/n}e^{i(1/n)\phi(z)},\quad z\in \Bbb D.\eqno(26)
     $$ 
Then  $F_n$ is a holomorphic function on $\Bbb D$  and  
 $\left(F_n\right)^n=F$ for all positive integers  $n$. 
Since  $F\in M^p\subset N^+$ and $F(z)\not=0$ for each  $z\in \Bbb D$, 
there exists $F^*(e^{i\theta})$ for almost 
every $e^{i\theta}\in\Bbb T$ and  $F^*(e^{i\theta})\not= 0$ 
for almost every $e^{i\theta}\in\Bbb T$.
For such a fixed  $e^{i\theta}$,  
 $\rho(re^{i\theta})$ is a positive  
continuous  function of  the variable  $r$ on the segment  $[0,1]$.
Therefore, there are  positive real numbers 
$l_{\theta}$ and $L_{\theta}$ depending on $\theta$ such that 
     $$
0<l_{\theta}\le \rho(re^{i\theta})\le L_{\theta}<\infty,\quad \mathrm{ for
\,\, each}\,\, r\in [0,1].\eqno(27)
     $$
Moreover, $\phi(re^{i\theta})$ is also a continuous  function of  
the variable  $r$ on the segment  $[0,1]$, 
 and hence, it is bounded on  $[0,1]$. This shows that
     $$
F_n(re^{i\theta})\to 1\quad{\mathrm as}\,\, n\to\infty,\eqno(28)
     $$  
uniformly on  $r\in [0,1]$. This together with the fact that
     $$
 Mf(\theta)=\sup_{0\le r<1}|f(re^{i\theta})|,\eqno(29)
     $$ 
yields
     $$
M(F_n-1)(\theta)\to 0\quad\mathrm{as}\,\, n\to\infty,\eqno(30)
     $$
for almost every  $e^{i\theta}\in\Bbb T$. Since $(F_n)^n=F$, we find that  
     $$
\log^+MF_n(\theta)\le\frac{1}{n}\log^+MF(\theta)\le \log^+MF(\theta),
\quad  n=1,2,\ldots ,\eqno(31)
     $$
Applying the  elementary  inequalities 
$M(F_n-1)(\theta)\le 1+ MF_n(\theta)$, $\log(2+a)\le 2\log 2+\log^+a$ $(a>0)$,
$(a+b)^p\le 2^{p-1}(a^p+b^p)$ ($a,b\ge 0, p>1$)
and the inequality (31), we find that
     \begin{eqnarray*}
\left(\log(1+M(F_n-1)(\theta))\right)^p &\le& \left(\log(2+MF_n(\theta))
\right)^p\\
\qquad\qquad\qquad  &\le& \left(2\log 2+\log^+MF_n(\theta)\right)^p\qquad\qquad 
\qquad\qquad\qquad\qquad(32)\\
&\le& 2^{p-1}\left((2\log 2)^p+(\log^+MF_n(\theta))^p\right)\qquad\qquad
\\
&\le& 2^{p-1}\left((2\log 2)^p+(\log^+MF(\theta))^p\right),
\quad n=1,2,\dots .
     \end{eqnarray*}     
By the inequality (32) and the fact that $F\in M^p$,
we conclude that  $F_n\in M^p$ for every positive integer $n$. 
Finally, by the  inequality (32) and  (30), we 
can apply the  Lebesgue dominated  convergence theorem to obtain
     $$
\int_0^{2\pi}\left(\log(1+M(F_n-1)(\theta))\right)^p\,\frac{d\theta}{2\pi}
\to 0\quad \mathrm{as}\,\, n\to\infty,\eqno(33)
     $$
or equivalently, $\rho_p(F_n,1)\to 0$ as $n\to\infty$. This completes 
the proof.\hfill $\square$
 \vspace{2mm}

The following result can be considered as the 
Hardy-Littlewood maximal theorem for holomorphic functions.

\vspace{2mm}

\noindent{\bf Lemma 9} (Hardy-Littlewood; see \cite[Theorem 1.9, p. 12]{d}). 
{\it Let  $f\in H^p$ with some $p$ such that $1< p<\infty$, and let
    $$
 Mf(\theta)=\sup_{0\le r<1}|f(re^{i\theta})|,\eqno(34)
   $$
Then $Mf\in L^p(\Bbb T)$ and 
    $$
\int_0^{2\pi}\left(Mf(\theta)\right)^p\,\frac{d\theta}{2\pi}\le 
C_p\int_0^{2\pi}|f^*(e^{i\theta})|^p\,\frac{d\theta}{2\pi},\eqno(35) 
   $$
where $C_p$  is a  constant  depending only on $p$.}

\vspace{2mm}

In order to prove Lemma 11, we will need the following known result.

\vspace{2mm}

\noindent{\bf Lemma 10} (see \cite[pp. 65--66, Proof of Lemma]{ho}).  
{\it Let  $B$ be an 
infinite Blaschke product with zeros $\{a_k\}_{k=1}^{\infty}$ 
and let   $\{B_n\}_{n=1}^{\infty}$ be a sequence of functions defined as 
      $$
  B_n(z)=\prod_{k=1}^{n}\frac{|a_k|}{a_k}\cdot\frac{a_k-z}{1-\bar{a}_kz},
\quad z\in\Bbb D. \eqno(36)
      $$
Then $B_n\to B$ in the space $H^2$, or equivalently,
      $$
 \int_0^{2\pi}|B(e^{i\theta})-B_n(e^{i\theta})|^2\,\frac{d\theta}{2\pi}
 \to 0\quad {\rm as}\quad n\to\infty.\eqno(37)
     $$

\vspace{2mm}

\noindent{\bf Lemma 11.} Let  $B$ be an infinite Blaschke product with zeros 
$\{a_k\}_{k=1}^{\infty}$. For each positive integer $n$,
take $B(z)=B_n(z)g_n(z)$ ($z\in\Bbb D$), where
      $$
  B_n(z)=\prod_{k=1}^{n}\frac{|a_k|}{a_k}\cdot\frac{a_k-z}{1-\bar{a}_kz},
\quad z\in\Bbb D\eqno(38)
      $$
and
       $$
g_n(z)=\prod_{k=n+1}^{\infty}\frac{|a_k|}{a_k}\cdot\frac{a_k-z}{1-\bar{a}_kz},
\quad z\in\Bbb D.\eqno(39)
       $$
 Then for each $p>1$,  $g_n\to 1$ in the space $M^p$ as $n\to\infty$.}

\vspace{2mm}

\noindent{\it Proof.} Using the inequality $\log(1+a)\le pa^{1/p}$  
for $a>0$ and the Cauchy-Schwarz integral inequality, we find that 
     \begin{eqnarray*}       
\rho_p(g_n,1)&=&\left(\int_0^{2\pi}\left(\log(1+M(g_n-1)(\theta))\right)^p\,
\frac{d\theta}{2\pi}\right)^{1/p}\\
\qquad\quad\qquad\quad\qquad\quad\qquad\quad &\le& p\left(\int_0^{2\pi}M(g_n-1)(\theta)\,\frac{d\theta}{2\pi}\right)^{1/p}
\qquad\quad\qquad\quad\qquad\quad (40) \\
&\le& p\left(\int_0^{2\pi}\left(M(g_n-1)(\theta)\right)^2\,
\frac{d\theta}{2\pi}\right)^{1/2p}.
     \end{eqnarray*}
Then using the inequality   (35) of Lemma 9, the well known fact 
that  $|B_n^*(e^{i\theta})|=1$ for almost every $e^{i\theta}\in\Bbb T$  and 
the inequality (40), we find that  
   \begin{eqnarray*}
\rho_p(g_n,1)&\le& p\left(C_2\right)^{1/{2p}}\left(\int_0^{2\pi}
|g_n(e^{i\theta})-1|^2\,\frac{d\theta}{2\pi}\right)^{1/{2p}}\\
\qquad\quad\qquad\quad\qquad &=& p\left(C_2\right)^{1/{2p}}\left(\int_0^{2\pi}
|B(e^{i\theta})-B_n(e^{i\theta})|^2\,\frac{d\theta}{2\pi}\right)^{1/{2p}}.
\qquad\quad\qquad\quad (41)
    \end{eqnarray*}
Finally, (37) of Lemma 10 and the inequality (41)
yield $\rho_p(g_n,1)\to 0$, or equivalently, 
$g_n\to 1$ in the space  $M^p$ as $n\to\infty$.  \hfill $\square$

\vspace{2mm}

\noindent{\it Proof of Theorem $5$.} We follow the proof of
Theorem 4 in \cite{k1}. Suppose that  ${\mathcal M}\not={\mathcal M}_{\lambda}$ for every 
$\lambda\in \Bbb D$. By Lemma 2, ${\mathcal M}$ contains 
a bounded holomorphic function $f$ on the disk $\Bbb D$. 
Then by Theorems 1 and 4,   $f$ can be factorized as a product $f=BF$, 
where  $B$ is the Blaschke product whose zeros are the same 
as these of $f$, and $F$ is a bounded holomorphic function
on $\Bbb D$ such that $F(z)\not= 0$ for all $z\in\Bbb D$.
Since  ${\mathcal M}$ is a prime ideal, we conclude that
either $F\in{\mathcal M}$ or $B\in{\mathcal M}$. Let $F\in {\mathcal M}$
and let $\{F_n\}_{n=1}^{\infty}$ be a sequence 
of functions defined as in Lemma 8.
Since $\left(F_n\right)^n=f$ for all $n=1,2,\ldots$, 
and taking into account that 
${\mathcal M}$ is a prime ideal, it follows that $F_n\in{\mathcal M}$ 
for all $n=1,2,\ldots$. This 
together with the fact that by Lemma 8, $F_n\to 1$ in 
$M^p$ as $n\to\infty$,   shows that  $1\in {\rm cl}({\mathcal M})$
(${\rm cl}({\mathcal M})$ is the closure of ${\mathcal M}$ in 
the space $M^p$). Hence, since ${\mathcal M}$ is also an ideal,
we conclude that  ${\rm cl}({\mathcal M})=M^p$, which is a contradiction 
in view 
of the assumption that  ${\mathcal M}$
is not dense in  $M^p$. This shows that it must be  $B\in{\mathcal M}$. Then 
we consider the cases when $B$ is a finite Blaschke product and when $B$ is 
an infinite Blaschke product. 

In the first case, set   
     $$
  B(z)=z^m\prod_{k=1}^{n}\frac{|a_k|}{a_k}\cdot\frac{a_k-z}{1-\bar{a}_kz}.\eqno(42)
     $$
Since $B\in{\mathcal M}$ and ${\mathcal M}$ is a prime ideal, it follows that 
either $z\in {\mathcal M}$ or $(a_k-z)/(1-\bar{a}_kz)\in{\mathcal M}$ 
for some   $k\in\{ 1,2,\ldots\}$. If $z\in {\mathcal M}$, then 
    $$
{\mathcal M}_0=zM^p\subset{\mathcal M}.\eqno(43)
    $$
Since by Theorem 4, $M^p=N^p$, it follows from \cite[Proposition 2.2]{msu}
that ${\mathcal M}_0$ is a maximal ideal in $M^p$. Hence, from (43) we 
conclude that ${\mathcal M}={\mathcal M}_0$. A contradiction. 

If  $(a_k-z)/(1-\bar{a}_kz)\in{\mathcal M}$ for
some   $k\in\{ 1,2,\ldots\}$, then for such a $k$ it must be
     $$
{\mathcal M}_{a_k}=\frac{a_k-z}{1-\bar{a}_kz}M^p\subset{\mathcal M}.\eqno(44)
    $$
Notice that in view of the fact that $M^p=N^p$,   
by \cite[Proposition 2.2]{msu}, the set 
$(a_k-z)M^p:=\{(a_k-z)f: f\in M^p \}$ 
is a  maximal ideal in $M^p$. This together with the fact that 
   $1/(1-\bar{a}_kz)\in H^\infty\subset M^p$ implies that  
${\mathcal M}_{a_k}$ is a maximal ideal in $M^p$. Therefore,  
the inclusion (44) yields ${\mathcal M}={\mathcal M}_{a_k}$. A contradiction.   

It remains to consider the case when $B$ is an infinite Blaschke product,
i.e.,
      $$
  B(z)=z^m\prod_{k=1}^{\infty}\frac{|a_k|}{a_k}\cdot\frac{a_k-z}{1-\bar{a}_kz}.\eqno(45)
     $$
 Define the sequence $\{g_n\}_{n=1}^{\infty}$ of functions  on $\Bbb D$ as 
   $$ 
 g_n(z)=\prod_{k=n+1}^\infty\frac{|a_k|}{a_k}\cdot\frac{a_k-z}{1-\bar{a}_kz},
z\in\Bbb D, \quad n=1,2,\ldots.\eqno(46)
    $$
Then  $g_n\in{\mathcal M}$ for all $n=1,2,\ldots$ because of the facts 
that  ${\mathcal M}$   is a prime ideal in $M^p$  and  $B\in M^p$. 
By Lemma 11,   $g_n\to 1$ in the space $M^p$ as $n\to\infty$.  
Therefore, it must be $1\in{\rm cl}({\mathcal M})$.
This  contradiction implies that   ${\mathcal M}={\mathcal M}_{\lambda}$
for some $\lambda\in \Bbb D$. This concludes proof of Theorem 5.
\hfill $\square$

\end{document}